\documentstyle[12pt]{article}
\def\baselinestretch{1.6}
\newtheorem{theorem}{Theorem}

\def\F{{\cal F}}

\def\P{{\cal P}}

\def\R{{\cal R}}

\def\endpf{\n$\Box$}

\def\be{\begin{equation}}
\def\ee{\end{equation}}

\countdef\eqn=101
\def\nr{\global\advance\eqn by 1
\eqno{(\the\eqn)}}

\newcommand{\btm}{\begin{itemize}}
\newcommand{\etm}{\end{itemize}}

\newcommand{\bc}{\begin{center}}
\newcommand{\ec}{\end{center}}
\newcommand{\ba}{\begin{array}}
\newcommand{\ea}{\end{array}}
\newcommand{\bey}{\begin{eqnarray}}
\newcommand{\eey}{\end{eqnarray}}
\newcommand{\ben}{\begin{eqnarray*}}
\newcommand{\een}{\end{eqnarray*}}
\newcommand{\bt}{\begin{tabular}}
\newcommand{\et}{\end{tabular}}
\newcommand{\n}{\noindent}
\def\Ref{\global\advance\refn by 1
\def\0{\par\hangindent 10pt\noindent\hskip 10pt\hskip -10pt}

\countdef\refn=112
\refn=0

\0 [\the \refn] }
\begin{document}
\bibliographystyle{plain}

\thispagestyle{empty}
\setcounter{page}{0}

{\Large L. Gy\"orfi, G. Morvai and  S. Yakowitz: }

\vspace {1cm}

{\Large Limits to consistent on-line forecasting for ergodic time series. }

\vspace {1cm}

{\Large IEEE Trans. Inform. Theory  44  (1998),  no. 2, 886--892.}

\vspace {2cm}

\begin{abstract}
This study concerns problems of time-series forecasting under
the weakest  of assumptions. Related results are surveyed
and are points of departure for the developments here,
some of which are new and others are new derivations
of previous findings. The contributions in this study are
all negative, showing that various plausible prediction
problems are unsolvable, or in other cases, are not solvable
by predictors which are known to be consistent when mixing
conditions hold.
\end{abstract}

\pagebreak


\section{Introduction}

\def\baselinestretch{1.6}

Given a 
 random variable sequence, such as  $X_0^{n-1}=(X_0,\dots,X_{n-1})$, a typical
prediction problem is to provide  from this data 
an estimate, say ${\hat E}(X_0^{n-1})$ of the succeeding value
$X_n.$
Following the influential book {\em Extrapolation, Interpolation,
and Smoothing of Stationary Time Series}
by N. Wiener \cite{Wiener}, the emphasis in prediction theory has
been (and still is) to find estimators which are convolutions

\be \hat E(X_0^{n-1})=\sum_{i=1}^{n} \alpha_i X_{n-i} \label{conv}\ee
of preceding observations.  Here the $\alpha_i\;'s$ are presumed to
be fixed real numbers determined entirely by the process covariance
function.    It is of course well-known that aside from the
Gaussian process case, linear predictors do not 
generally give the least-squares optimal prediction, or even a consistent
approximation (as the data base grows) of the optimal estimator, which
is the conditional expectation $E(X_n|X_0^{n-1})$ of $X_n$. 
If the time series happens to be generated by the
nonlinear autoregression $X_{n}=\sqrt{|X_{n-1}|}+\epsilon_n$
for some i.i.d. non-singular noise sequence $\{\epsilon_n\}$, then
no matter how the linear parameters in (\ref{conv}) are adjusted,
the expected squared-error prediction of $X_{n}|\{X_i, i < n\}$
 will be worse than the estimate
$m(X_{n-1})=\sqrt{|X_{n-1}|}.$

The Kalman filter and ARMA (or as it is sometimes called,
 Box/Jenkins) methods are equivalent to (\ref{conv}), as are
predictors based on spectral analysis.
These ''second-order" techniques were well-suited to the 
period before about 1970 when
data set size and access to computer power were relatively limited.

Beginning with  the pioneering work of 
Roussas \cite{Rous} and Rosenblatt \cite{Rosen}, nonparametric
methods worked their way into the literature of forecasting for dependent
series. Several people, including the authors,
 have  investigated  forecasting
problems, such as enunciated by Cover
\cite{Cover75}, under the sole hypotheses of stationarity and ergodicity.
Two classical results for stationary ergodic sequences,
namely,  Birkhoff's Theorem,
\[
\lim_{n\to\infty}{1\over n}\sum_{i=1}^n X_i=E(X) \ \ 
\mbox{almost surely,}
\]
and the Glivenko-Cantelli Theorem,

\[
\lim_{n\to\infty}\sup_x|F_n(x)-F(x)|=0  \ \ 
\mbox{almost surely,}
\]
 for convergence of the empirical to the true distribution function
are clear evidence that some statistical problems are 
solvable under weak assumptions regarding dependency.  
In fact, since nonergodic stationary sequences can be viewed as
mixtures of ergodic modes, ergodicity itself is not
a vital assumption for prediction.  This matter is discussed in \cite{MoYaAl96}.

 On the other hand, not all problems solvable for independent
sequences can be mastered in the general setting.  For instance,
 Gy\"orfi and Lugosi \cite{GyLu92} show that the kernel density
estimator is not universally consistent, even though we do have
consistency of the recursive kernel density estimator
under ergodicity provided that for some integer
$m_0$ the conditional density of 
$X_{m_0}$ given the condition $X^{0}_{-\infty}$ exists (Gy\"orfi and Masry \cite{GyMa90}).

It will be useful to distinguish between two classes of prediction problems.\\

\noindent{\bf Static forecasting.} Find an estimator
$\hat E(X_{-n}^{-1})$ of the value $E(X_0|X_{-N}^{-1})$ such that 
for  any stationary and ergodic sequence $\{X_i\}$ with values
in some given coordinate set $\cal X,$ 
almost surely, 
\be 
\label{ascvexp2}
\lim_{n\to\infty} \hat E(X_{-n}^{-1})=E(X_0|X_{-N}^{-1}). 
\ee
In (\ref{ascvexp2}), $N$ may be $\infty,$ in which case
we will speak of the  {\em static total-past} prediction. Otherwise,
this is called the  {\em static autoregression} problem.
In either case, it is presumed that
the forecaster $\hat E(X_{-n}^{-1}) $ depends only on the data
segment $X_{-n}^{-1}.$\\

\noindent
The other   problem of interest here is,\\

\noindent
{\bf Dynamic forecasting.} Find an estimator
$\hat E(X_{0}^{n-1})$ of the value $E(X_n|X_{n-N}^{n-1})$ such that 
for  any  stationary and ergodic sequence $\{X_i\}$ taking values
in a given set $\cal X,$
almost surely, 
\be 
\label{ascvexp1}
\lim_{n\to\infty} |\hat E(X_{0}^{n-1})-E(X_n|X_{n-N}^{n-1})|=0. 
\ee
Here $N$ is typically either $n$ or  a fixed postive
integer, and 
the estimator must be constructible from data collected
from time $0$ up to the ''current" time $n-1.$ When $N$ is  a fixed postive
integer,
we have the  {\em dynamic autoregression problem,} and
the alternative  category will be referred to as   the {\em dynamic total-past}
forecasting  problem.\\

When the coordinate set $\cal X$ is finite or countably infinite, for 
both autoregression problems ($N<\infty$)
 one may construct
an estimator  with consistency verified by simple application
of the ergodic theorem. Thus, for static autoregression,
 the observed sequence $X_{-N}^{-1}$ 
has positive marginal probability.  Define for $ n > N,$

\bey
Num(X_{-N}^{-1},n)&=&\sum_{j=1}^{n-N}  I_{[X^{-j-1}_{-j-N}=X_{-N}^{-1}]}\,X_{-j}\\
Denom(X_{-N}^{-1},n)&=&\sum_{j=1}^{n-N}  I_{[X^{-j-1}_{-j-N}=X_{-N}^{-1}]}\\
g(x_{-1},\dots,x_{-N})&=&E(X_0 1_{[X_{-N}^{-1}=x_{-N}^{-1}]})\\
h(x_{-1},\dots,x_{-N})&=&P(X_{-N}^{-1}=x_{-N}^{-1}).
\eey
From the ergodic theorem, a.s., 
\bey
{1\over n-N} Num(X_{-N}^{-1},n) \to g(X_{-N}^{-1})\\
{1\over n-N} Denom(X_{-N}^{-1},n) \to h(X_{-N}^{-1})
\eey
and this implies the consistency of the estimate 
$$\hat E(X_{-n}^{-1})=Num(X_{-N}^{-1},n)/Denom(X_{-N}^{-1},n),$$ i.e.,
almost surely, as $n\to\infty,$

\be \hat E(X_{-n}^{-1})\to 
{g(X_{-N}^{-1})\over h(X_{-N}^{-1})}=
E(X_0|X_{-N}^{-1}).\ee

For the dynamic case, take

\be
\hat E(X_0^{n-1})=\frac{\sum_{j=N}^{n-1} I_{\{X_{j-N}^{j-1}=X_{n-N}^{n-1}\}}X_{j}}
{\sum_{j=N}^{n-1} I_{\{X_{j-N}^{j-1}=X_{n-N}^{n-1}\}}}\ee
Since now there are but finitely many possible strings
$X_{n-N}^{n-1}$, the ergodic theorem implies we have
a.s. convergence of the estimator of the successor value on each of
them.

In  1978, Ornstein \cite{Ornstein78} provided
an estimator for the static, finite $\cal X$ total-past prediction problem.
In 1992, Algoet \cite{Algoet92} generalized
Ornstein's findings to
allow that $\cal X$ can be any Polish space.
  More recently, Morvai, Yakowitz and Gy\"{o}rfi  \cite{MoYaGy95} gave a
simpler algorithm and convergence proof for that problem.
 It is to be
admitted that
at this point,  these algorithms are terribly unwieldy.

The {\it partitioning} estimator is a representative computationally feasible nonparametric algorithm.  
Such methods attracted a great deal of theoretical attention in 
the 1980's, much of it being summarized and referenced in
the monograph \cite{GyHaSaVi89}.  This partitioning  method, 
and its relatives
such as the nearest neighbor and the kernel autoregressions,
are  known to  consistently estimate the conditional expectation
$E(X_0|X_{-1})$
under a great many ``mixing" conditions regarding the
degree of dependency  of the present and future
on the distant past cf. Chapter III. in \cite{GyHaSaVi89}.  These  mixing conditions, while plausible,
are difficult to check.  There is virtually no literature
on inference of mixing conditions and mixing parameters from data.

In view of these positive results under mixing, we wanted to show that the
partitioning regression estimate, known to be
effective for time series under a variety of mixing conditions, 
 suffices for static autoregressive forecasting, when
${\cal X}$ is real. Such a finding would be 
interesting because this method is straightforward to apply
and in a certain sense, is economical with data.
This conjecture turns out to be
untrue.  We will show that there exists a partition sequence which satisfies the usual conditions and a stationary ergodic time series $X_n$ such that   on a set of positive probability, for the partitioning estimate $\hat E(X_{-n}^{-1})$, 
\be 
\label{ascvexp4}
\limsup_{n\to\infty} |\hat E(X_{-n}^{-1})-E(X_0|X_{-1})|>0.
\ee
This and a related result are demonstrated in Section \ref{secpartition}.

Turning attention to  dynamic forecasting, in Section \ref{secdynamic}, we
relate a theorem due to Bailey \cite{Bailey76} stating that, 
in contrast to the
static case, even for binary sequences,
there is no algorithm that can achieve a.s. convergence
in the sense of (\ref{ascvexp1}), for the dynamic total-past problem with $N=n$.  
On the other hand,
it is evident that algorithms such as \cite{Algoet92}
or \cite{MoYaGy95}, which provide solution to the 
a.s. static forecasting problem can be modified to achieve
convergence in probability for this recalcitrant case.
Details of a conversion were given in \cite{MoYaAl96},
which gives yet another plan for attaining weak convergence
of dynamic forecasters.  
When the coordinate space is finite, it turns out 
that implicitly,  algorithms
for inferring entropy (e.g., \cite{ZiLe78}) can also be  utilized for
constructing weakly convergent static and dynamic autoregressive forecasters.
This has  been noted (e.g., \cite{Ryabko88}), and discussed at length 
in Section IV of \cite{MoYaAl96}.

\section {Dynamic forecasting} \label{secdynamic}

Let $\{X_i\}_{-\infty}^{\infty}$ be a stationary ergodic 
binary-valued process. The goal is to find a predictor 
$\hat E(X_{0}^{n-1})$ of the value $E(X_n|X_0^{n-1})$ such that 
almost surely, 
\[ 
\lim_{n\to\infty} |\hat E(X_{0}^{n-1})-E(X_n|X_0^{n-1})|=0
\]
for all stationary and ergodic processes. We show by
the statement below that this goal is not achieveable.

\noindent
\begin{theorem} \label{Th1}
{\sc (Bailey \cite{Bailey76}, Ryabko \cite{Ryabko88})}
 For any estimator 
$\{\hat E(X_{0}^{n-1})\}$ there is a stationary ergodic  binary-valued
process $\{X_i\}$ such that
\[
P(\limsup_{n\to\infty} 
|\hat E(X_0^{n-1})- E(X_n|X_0^{n-1})|\ge 1/4) 
\ge {1\over 8}.
\]
\end{theorem}
{\sc Remark}
Bailey's counterexample for dynamic total-past forecasting
 uses the technique of
cutting and stacking developed by  Ornstein \cite{Ornstein74}
(see also Shields \cite{Shields91}). Bailey's  proof has not been published 
and is hard to follow, whereas  Ryabko  omitted his lengthy proof and only 
sketched an intuitive argument in his paper. 
These results are not widely known. 
In view of their significance to the issue of the 
"limits of forecasting", we wanted to unambigously enter it into the 
easily-accessible literature. 

\noindent{\sc Proof} The present  proof is a simplification of the clever counterexample of Ryabko
 \cite{Ryabko88}.  First we define a 
Markov process which serves as the technical tool for  construction of our 
counterexample. Let the state space $S$ be the non-negative integers. From state $0$ the process certainly passes to state $1$ and then 
to state $2$, at the following epoch. From each state $s\ge 2$,  the 
Markov chain  passes either to state $0$ or to state 
$s+1$ with equal probabilities $0.5$.  
This construction yields a stationary and ergodic 
Markov process $\{M_i\}$ with stationary distribution
$$
P(M_i=0)=P(M_i=1)={1\over 4}
$$
and
$$
P(M_i=j)={1\over  2^{j}} \mbox{\ \  for $j\ge 2$}.
$$

\noindent
Let $\tau_k$ denote the first positive time of  occurence of 
state $2k$ :  
$$ \tau_k=
\min\{i\ge 0: M_i=2k\}.
$$
Note that if $M_0=0$ then $M_i\le 2k$ for $0\le i\le \tau_k$. 
Now we define the hidden Markov chain   $\{X_i\}$, which we denote as, $X_i=f(M_i)$. It will serve as the stationary unpredictable time series. 
We will use the notation $M_0^n$ to denote the
 sequence of states $M_0,\dots,M_n$. 
Let $f(0)=0$, $f(1)=0$, and $f(s)=1$ for all even states $s$. 
A feature of this definition of $f(\cdot)$ is that whenever  $X_n=0,X_{n+1}=0,X_{n+2}=1$ we  know 
 that $M_n=0$ and {\it vice versa}. 
Next we will define $f(s)$ for odd states $s$ maliciously. 
We define $f(2k+1)$ inductively for $k\ge 1$. Assume $f(2l+1)$ is defined for $l<k$. If $M_0=0$ (that is, $f(M_0)=0$, $f(M_1)=0$, $f(M_2)=1$) then 
$M_i\le 2k$ for $0\le i\le \tau_k$ and 
the mapping 
$$M_0^{\tau_k}\rightarrow(f(M_0),\dots,f(M_{\tau_k}))$$
is invertible. ( Given $X_0^n$ find $1\le l\le n$,  and  positive integers $0=r_0<r_1<\dots<r_l=n+1$ such that 
$X_0^n=(X_{r_0}^{r_1-1}, X_{r_1}^{r_2-1},\dots,X_{r_{l-1}}^{r_l-1})$, where 
 $2\le r_{i+1}-1-r_i<2k$ for $0\le i<l-1$, $r_l-1-r_{l-1}=2k$  and  for $0\le i<l$,  $X_{r_i}^{r_{i+1}-1}=(f(0),f(1),\dots,f(r_{i+1}-1-r_i))$.  
Now  $\tau_k=n$ and $M_{r_i}^{r_{i+1}-1}=(0,1,\dots,r_{i+1}-1-r_i)$ for $0\le i< l$.  This construction is always possible under our postulates that  $M_0=0$ and  $\tau_k=n$.)
Let 
$$
B_k^+= \{M_0=0, \hat E(f(M_0),\dots,f(M_{\tau_k}))\ge {1\over 4} \}
$$
and 
$$
B_k^-= \{M_0=0, \hat E(f(M_0),\dots,f(M_{\tau_k}))<  {1\over 4} \}.
$$
 Now notice that the events $B_k^+$ and $B_k^-$ 
do not depend on the future values of $f(2r+1)$ for $r\ge k$, and    
one of  these  events must have probability at least $1/8$ 
since $$
P(B_k^+)+P(B_k^-)=P(M_0=0)={1\over 4}.
$$
Let $I_k$ denote the most likely of the events $B_k^+$ and $B_k^-$, and 
inductively define 
\[
f(2k+1)= \left\{ \begin{array}{ll} 1 & \ \ \mbox{if $I_k=B_k^-$,} \\
0 & \ \ \mbox{if $I_k=B_k^+$.}
\end{array} \right. 
\]
  
\noindent 
  Because of the construction
of $\{M_i\}$,  on event $I_k$,
\begin{eqnarray*}
E(X_{\tau_k+1}|X_0^{\tau_k})&=&
f(2k+1) P(X_{\tau_k+1}=f(2k+1)| X_0^{\tau_k})\\
&=& f(2k+1) P(M_{\tau_k+1}=2k+1| M_0^{\tau_k})\\
&=& 0.5 f(2k+1).
\end{eqnarray*}
The conditional expectation $E(X_{\tau_k+1}|X_0^{\tau_k})$  and the 
estimate $\hat E(X_0^{\tau_k})$ differ at least $1/4$ on the event $I_k$ 
and this event occurs with probability at least $1/8$. By Fatou's lemma, 
\begin{eqnarray}
\nonumber
\lefteqn{P(\limsup_{n\to\infty} 
\{|\hat E(X_0^{n-1})- E(X_{n}|X_0^{n-1})|\ge 1/4\})}\\
\nonumber
&\ge&
P(\limsup_{n\to\infty} 
\{|\hat E(X_0^{n-1})- E(X_{n}|X_0^{n-1})|\ge 1/4,X_0=X_1=0, X_2=1\})\\
\nonumber
&\ge& P(\limsup_{k\to\infty} 
\{|\hat E(f(M_0),\dots,f(M_{\tau_k}))- E(f(M_{\tau_k+1})|f(M_0),\dots,f(M_{\tau_k}))|\ge 1/4, M_0=0\})\\
\nonumber
&\ge&
P(\limsup_{k\to\infty} I_k)= 
E(\limsup_{k\to\infty} 1(I_k))
\ge \limsup_{k\to\infty} E1(I_k)= \limsup_{k\to\infty} P(I_k)\ge {1\over 8}.
\end{eqnarray} 
\endpf

\bigskip

 We noted in the Introduction that 
there are static total-past  empirical forecasters
(i.e., $N=\infty$ in (\ref{ascvexp2})) which are 
strongly universally consistent  when the coordinate
space $\cal X$ is real.  These are readily transcribed
to weakly-consistent dynamic forecasters.  The following
(which was  inspired by the methods of \cite{Ryabko88})
shows that one cannot hope for a strongly consistent 
autoregressive dynamic forecaster.

Let $\{X_i\}_{-\infty}^{\infty}$ be a stationary ergodic 
real-valued process. The goal is to find a one-step predictor 
$\hat E(X_{0}^{n-1})$ of the value $E(X_n|X_{n-1})$ (i.e. $N=1$) such that 
almost surely, 
\[ 
\lim_{n\to\infty} |\hat E(X_{0}^{n-1})-E(X_n|X_{n-1})|=0
\]
for all stationary and ergodic processes. 

\begin{theorem} {\sc (Ryabko \cite{Ryabko88})} \label{Th2}
 For any estimator 
$\{\hat E(X_{0}^{n-1})\}$ there is a stationary ergodic process $\{X_i\}$ with values from a countable subset of the real numbers  such that
\[
P(\limsup_{n\to\infty} 
\{|\hat E(X_0^{n-1})- E(X_n|X_{n-1})|\ge 1/8\}) 
\ge {1\over 8}.
\]
\end{theorem}
{\sc Proof}
We will use the  Markov process $\{M_i\}$
defined in the proof of Theorem \ref{Th1}.
Note that one must pass through state $s$ to get to any
state $s'>s$ from $0$.  
We construct a hidden Markov chain $\{X_i\}$ which is 
in fact just a relabeled version of $\{M_i\}$.
This construct uses  a different (invertible) 
function $f(\cdot),$ for $X_i=f(M_i)$.
Define f(0)=0, $f(s)=L_s+2^{-s}$ if $s>0$ where 
$L_s$ is either $0$ or $1$ as specified later. In this way, knowing $X_i$ is equivalent to knowing $M_i$ 
and {\em vice versa}.  Thus $X_i=f(M_i)$ where $f$ is one-to-one. 
For $s\ge 2$ the conditional expectation is,
$$E(X_t|X_{t-1}=L_s+2^{-s})={L_{s+1} + 2^{-(s+1)}\over 2}.$$ 
We complete the description of 
the function $f(\cdot)$ and thus the
conditional expectation by defining
 $L_{s+1}$ so as to confound any proposed 
predictor $\hat E(X_0^{n-1}).$
Let  $\tau_s$ denote the time of first occurence of state $s$ :  
$$ \tau_s=
\min\{i\ge 0: M_i=s\}
$$
Let $L_1=L_2=0$. Suppose $s\ge 2$. Assume we specified $L_i$ for $i\le s$. Define  
$$
B_s^+=\{X_0=0, \hat E(X_0^{\tau_s})\ge {1\over 4} \}
$$
and 
$$
B_s^-=\{X_0=0, \hat E(X_0^{\tau_s})< {1\over 4} \}.
$$
One of the two events must have probability at least $1/8$.  
Take $L_{s+1}=1$, and $I_s=B_s^-$  if $P(B_s^-)\ge P(B_s^+).$ 
Let $L_{s+1}=0$, and $I_s=B_s^+$ if  $P(B_s^-)< P(B_s^+)$. 
The difference of the estimate and the conditional expectation is at 
least $1/8$ on the  event  $I_s$ and this event occurs  with probability  
not less  than  $1/8$. 
By Fatou's lemma, 
\begin{eqnarray}
\nonumber
\lefteqn{P(\limsup_{n\to\infty} 
\{|\hat E(X_0^{n-1})- E(X_n|X_{n-1})|\ge {1\over 8}\})}\\ 
\nonumber
&\ge&
P(\limsup_{s\to\infty} 
\{|\hat E(X_0^{\tau_s})- E(X_{\tau_s+1}|X_{\tau_s})|\ge {1\over 8},X_0=0\})\\ 
\nonumber
&\ge&
P(\limsup_{s\to\infty} I_s) \ge \limsup_{s\to\infty} P(I_s) \ge {1\over 8}.
\end{eqnarray} 
\endpf

\noindent
{\bf Remark 1.} The counterexample  in Theorem~\ref{Th2} is a 
Markov chain with countable number of states. 
(The correspondence between states $s$ and labels $f(s)$ is one-to-one.)

\medskip

\noindent
{\bf Remark 2.} One of the referees noted  that the question of whether 
strongly consistent forecasters exist if the process is postulated to be 
Gaussian, is interesting and open.    

\section{ Partitioning estimates which are  not universally consistent for autoregressive
static forecasting} \label{secpartition}

Let $\{(Y_i,Z_i)\}_{-\infty}^{\infty}$ be a  
stationary sequence taking values from $\R\times\R$. 
Let $\P_n=\{A_{n,j}\}$ be a partition of the real line. Let $A_n(z)$ denote
the cell $A_{n,j}$ of $\P_n$ into which $z$ falls. 
Let 
 \be
\nu_n(A)={1\over n-1} \sum_{i=1}^{n-1} I_{[Z_{-i}\in A]}Y_{-i}
\ee
and 
\be
\mu_n(A)={1\over n-1} \sum_{i=1}^{n-1} I_{[Z_{-i}\in A]}. 
\ee
Then the {\it partitioning estimate} of the regression function $E(Y_0|Z_0=z)$ 
is defined as follows:
\be 
\label{defpartesthmnz}
\hat m_n(z)={\nu_n(A_n(z))\over \mu_n(A_n(z))}=
 { \sum_{i=1}^{n-1} I_{[Z_{-i}\in A_n(z)]} Y_{-i} \over  
\sum_{i=1}^{n-1} I_{[Z_{-i}\in A_n(z)]}}.
\ee
We follow the convention that $0/0=0$. 

If $\{(Y_i,Z_i)\}$ is i.i.d.  or uniform
mixing or strong mixing with  certain assumptions on the
rates of the mixing parameters,
 then the strong universal consistency of
the partitioning estimate has been demonstrated
under the proviso that for all intervals $S$ symmetric 
around $0$,
\be
\label{r1}
\lim_{n\to \infty} \sup_{j;A_{n,j}\cap S\ne \phi}diam(A_{n,j})=0
\ee
and
\be
\label{r2}
\lim_{n\to \infty} {|\{j;A_{n,j}\cap S\ne \phi \}|\over n}=0
\ee
(cf. Devroye  and Gy\"orfi \cite{DeGy83} and Gy\"orfi \cite{Gyor91},
for the i.i.d. case, and Chapter III. in \cite{GyHaSaVi89} for mixing and for cubic partitions).

In the discussion to follow, we investigate the problem of one-step (i.e. $N=1$) 
autoregressive static
forecasting by  the partitioning estimate 
for the case of a stationary and ergodic 
real-valued 
process $\{X_i\}_{-\infty}^{\infty}$.  Thus the intention is to 
infer the value 
$m(x)=E(X_0|X_{-1}=x)$. In this case the  partitioning estimate 
is adapted for autoregressive prediction.  The predictor 
$\hat m_n(x)$ 
is here defined to be the partitioning estimate 
$\hat m_n(z)$ in (\ref{defpartesthmnz}) with  $z=x$ for the process 
$\{Y_i=X_i,Z_i=X_{i-1}\}_{-\infty}^{\infty}$. That is, 
 \be
\hat m_n(x)={\nu_n(A_n(x))\over \mu_n(A_n(x))}=
{\sum_{i=1}^{n-1} I_{[X_{-1-i}\in A_n(x)]}X_{-i}\over 
 \sum_{i=1}^{n-1} I_{[X_{-1-i}\in A_n(x)]}}.
\label{partit}
\ee
In an obvious way, the partitioning estimate results in a one-step static forecasting: $\hat E(X_{-n}^{-1})=\hat m_n(X_{-1})$. 

In contrast to the success of the partitioning estimate 
for independent   or   mixing sequences, we have the following negative results. 

\noindent
\begin{theorem} \label{Th3}

There is a stationary ergodic process $\{X_i\}$ with marginal distribution
uniform on $[0,1)$
and a sequence of partitions $\P_n$ 
satisfying (\ref{r1}) and (\ref{r2}) such that for the partitioning
forecaster $\hat m_n(X_{-1}),$  defined by (\ref{partit}),
$$
P(\limsup_{n\to\infty} |\hat m_n(X_{-1})-m(X_{-1})| \ge 0.5)\ge 0.5. 
$$
\end{theorem}
{\sc Proof}
We will construct a sequence of subsets  $B_n$ of $[0,1)$,  such that
$$P(X_{-1}\in \limsup_{n\to\infty} B_n)>0$$ and  if 
$X_{-1}\in B_n$ then $X_{-2}\notin B_n, \dots, X_{-n}\notin B_n$.  
Thus, when $X_{-1}\in B_n,$ we will be assured that none of
the data values up to time $n$ are in this set, and consequently a 
conventional partitioning estimate has no data in the appropriate
partition cell. 
We present first a dynamical system. We will define a transformation $T$ on the unit interval.
 Consider the binary expansion 
$r_1^{\infty}$  of each real-number $r\in [0,1)$, that is, 
$r=\sum_{i=1}^{\infty} r_i 2^{-i}$.  When there are two expansions,
use the representation which contains finitely many  $1's$.
Now let 
\be
\tau(r)= \min\{i>0: r_i=1\}.
\ee
Notice that, aside from the exceptional set  $\{0\}$,  which has Lebesgue measure zero $\tau$  is finite and 
well-defined on the closed unit interval.   
The transformation is defined by 
\be
(Tr)_i=\left\{
\begin{array}{ll}
1 & \mbox{if $0<i<\tau(r)$} \\
0 & \mbox{if $i=\tau(r)$} \\
r_i & \mbox{if $i>\tau(r)$}.
\end{array}
\right. 
\ee
Notice that in fact, $Tr=r-2^{-\tau(r)}+\sum_{l=1}^{\tau(r)-1} 2^{-l}$.
  All iterations $T^k$ of $T$ for $-\infty<k<\infty$ are well defined and invertible with the exeption of the set of dyadic rationals which has Lebesgue measure zero. In the future we will neglect this set.  
One of the referees pointed out that transformation $T$ could be  defined recursively as 
$$
Tr=\left\{
\begin{array}{ll}
r-0.5 & \mbox{if $0.5\le r<1$} \\
{1+T(2r) \over 2} & \mbox{if $0\le r< 0.5$.} 
\end{array}
\right.
$$
Let $S_i=\{
I_0^i,\dots,I_{2^i-1}^i\}$ be a partition of $[0,1)$ 
where for each integer $j$ in the range $0\le j<2^i$ $I_j^i$ is defined as the set of numbers $r=\sum_{v=1}^{\infty} r_v 2^{-v}$ whose 
binary expansion $0.r_1, r_2,\dots$ starts with the bit sequence $j_1,j_2,\dots,j_i$ that is reversing the binary expansion $j_i,\dots,j_2,j_1$ of the number 
 $j=\sum_{l=1}^i 2^{l-1} j_l$. 
  Observe that in $S_i$ there are $2^i$ left-semiclosed 
intervals and each interval $I_j^i$ has length (Lebesgue measure) $2^{-i}$. 
 Now $I_{j}^i$ is mapped linearly, under $T$ onto 
$I_{j-1}^i$ for $j=1,\dots,2^i-1$. 
To confirm this, observe that  for $ j=1,\dots,2^i-1$, if $r\in I_{j}^i$ then 
\begin{eqnarray*}
Tr&=&\sum_{l=1}^{\tau(r)-1} 2^{-l}+\sum_{l=\tau(r)+1}^{\infty} r_l 2^{-l} \\
&=&r- \sum_{l=1}^i 2^{-l} ( j_l-(j-1)_l)\\
&=& \sum_{l=1}^{i} (j-1)_l 2^{-l}+\sum_{l=i+1}^{\infty} r_l 2^{-l}.  
\end{eqnarray*}
Now if $0<r\in I_0^i$ then $\tau(r)>i$ and  so $Tr\in I_{2^i-1}^i$. 
Furthermore, if $r\in I_{2^i-1}^i$ then 
$r_1=\dots=r_i=1$,  and thus conclude that $(T^{-1}r)_1=\dots=(T^{-1}r)_i=0$, that is, $T^{-1}r\in I_0^i$. Let $r\in [0,1)$ and $n\ge 1$ be arbitrary. Then   $r\in I_j^n$ for some $0\le j\le 2^n-1$. 
For all $j-(2^n-1)\le k\le j$,  
 \begin{equation} \label{iteratedshift}
T^{k}r= \sum_{l=1}^{n} (j-k)_l 2^{-l}+\sum_{l=n+1}^{\infty} r_l 2^{-l}.  
\end{equation} 
Now since $T^{-1}I_{j}^i=I_{j+1}^i$ 
for $i\ge 1$, $j=0,\dots, 2^i-2$, and the union over $i$ and $j$ of these 
sets generate the Borel $\sigma$-algebra, 
we conclude that $T$ is measurable. Similar reasoning shows that $T^{-1}$ is also measurable.  
The dynamical system $(\Omega,\F,\mu,T)$ is identified with
$\Omega=[0,1)$ and $\F$  the Borel $\sigma$-algebra on $[0,1)$, 
 $T$ being the 
transformation developed above.  Take $\mu$ to be 
 Lebesgue measure on the 
unit interval. Since transformation $T$ is measure-preserving on each set 
in the collection  $\{I_j^i : 1\le j\le 2^i-1, 1\le i<\infty\}$ 
and these intervals generate the Borel $\sigma$-algebra $\F$, 
$T$ is a stationary transformation. Now we prove that  transformation $T$ is ergodic as well. 
Assume $TA=A$. If $r\in A$ then $T^l r\in A$ for $-\infty<l<\infty$. Let $R_n: [0,1)\rightarrow \{0,1\}$ be the function $R_n(r)=r_n$. If $r$ is chosen uniformly on $[0,1)$ then $R_1,R_2,\dots$ is a series if i.i.d. random variables. 
Let $\F_n=\sigma(R_n,R_{n+1},\dots)$. By~(\ref{iteratedshift}) it is immediate that $A\in \cap_{n=1}^{\infty} \F_n$ and so  $A$ is a tail event. By  Kolmogorov's zero one law $\mu(A)$ is either zero or one.  Hence $T$ is ergodic. 

Next we construct the sequence $\{B_n\},$ described at the beginning
of this proof, which 
forces the partitioning method to make``no data''  estimations infinitely
often. For each $B_n$ we require that 
\begin{equation}
\label{condBndisjoint}
T^0 B_n,\dots,T^{-n} B_n \ \ \mbox{ be disjoint.}
\end{equation}
The definition is inductive on  $k\ge 1$. 
For   $k=1,$ we define  $B_1=I_0^1$, that is $B_1$ is taken to be the left 
half of the unit interval. 
Since $T^{-1} I_0^1=I_1^1$ condition (\ref{condBndisjoint}) is satisfied. 
Recursively, for $k=2,3,\ldots$ we define $B_l$  
for $2^{k-2}< l\le 2^{k-1}$. 
Suppose that by the end of the construct for $k-1$ we have defined 
$B_l$ for 
$1\le l\le 2^{k-2}$ so that condition (\ref{condBndisjoint}) is satisfied 
with $n=l$. 
For the next iteration, $k,$ we define $B_{2^{k-2}+l}$ for $1\le l\le 2^{k-2}$ by
$$B_{2^{k-2}+l}=I_{2^{k-1}-2l}^k$$ 
and since 
$$
T^{-m} B_{2^{k-2}+l}=I_{2^{k-1}-2l+m}^k
$$
for $0\le m\le 2^{k-2}+l$, condition (\ref{condBndisjoint}) is satisfied.
Take $C_k$ to be the union of the newly defined $B_l's$:
$$C_k=\bigcup_{2^{k-2}<l\le 2^{k-1}} B_l=
\{r=0.r_1,\dots : r_1=0, r_k=0\}.
$$
Now 
\begin{eqnarray*}
\mu(\limsup_{n\to\infty} B_n)&=&
\mu(\limsup_{n\to\infty} C_n)\\
&=&\mu(\{r\in [0,1) : r_1=0, r_n=0 \ \mbox{for infinitely many n}\})\\
&=& 0.5
\end{eqnarray*}
since the set of real numbers in $[0,0.5)$ having infinitely many zero bits 
in their expansion constitute a set of Lebesgue measure $0.5$. 
Define the process as follows: For $\omega$ randomly chosen
from $[0,1)$ according to Lebesgue measure $\mu,$ the
dynamical system construct has us take,
$X_i(\omega)=T^{i+1}\omega$. 
Notice that the time series $\{X_i\}_{-\infty}^{\infty}$ is not just 
stationary and ergodic but also Markovian with continuous state space. 
Notice also that any observation $X_i$ determines 
the entire future and  past. 
By (\ref{condBndisjoint}) if 
$\omega\in B_n$ then $X_{-1}(\omega)\in B_n$ and 
$X_{-i}(\omega)\notin B_n$ for all $1<i\le n$.  
We will construct  a  partitioning estimator which satisfies the conditions
of the definition given above and yet which is ineffective for this process. 
Take  $\{H_{n,j}\}_{j=1}^{q(n)}$ to be a partition of $[0,1)$ by intervals 
of length $h_n=1/q(n)$ such that 
\begin{equation}
\label{condhntozero}
h_n\to 0
\end{equation}
and
\begin{equation}
\label{condnhntoinfty}
n h_n \to \infty.
\end{equation}  
Let $A_{n,j}^+=H_{n,j}\cap B_n$ and 
$A_{n,j}^-=H_{n,j}\cap \bar{B_n}$,  the overbar denoting
complementation.
Choose $\P_n=\{A_{n,j}^+,A_{n,j}^- \mbox{\  : \ $j=1,\dots,q(n)$}\}$. 
Partition $\P_n$ satisfies the conditions  (\ref{r1}) and (\ref{r2}). 
If $\omega\in B_n$ then for some $1\le j\le q(n)$, 
$X_{-1}(\omega)\in A_{n,j}^+$ and 
$X_{-i}(\omega)\notin A_{n,j}^+$ for all $1<i\le n$.  
The left half $B_1=I_0^1$ of $[0,1)$ is mapped to the right half $T B_1=I_1^1$
and $B_n\subseteq B_1$, so $E(X_0|X_{-1})(\omega)\ge 0.5$ if  $\omega\in B_n$. 
On the other hand,  $\hat m_n(X_{-1})(\omega)=0$ if 
$\omega\in B_n$. Thus 
$$
P(\limsup_{n\to\infty} |\hat m_n(X_{-1})-m(X_{-1})| \ge 0.5)
\ge \mu (\limsup_{n\to\infty} B_n)=0.5.
$$
\endpf

\bigskip

\begin{theorem} \label{Th4} For the partitioning estimate 
 $\hat m_n(x)$, defined by (\ref{partit}), there is a stationary ergodic process $\{X_i\}$ with marginal distribution
uniform on $[0,1)$
and a sequence of partitions $\P_n$ 
satisfying (\ref{r1}) and (\ref{r2}) such that for large  $n$, 
\[
P(\int |\hat m_n(x)-m(x)|\mu(dx)\ge 1/16)\ge {1\over 8}. 
\]
\end{theorem}

\noindent
{\sc Proof } The proof is a slight extension of the Shields' construction where he proved the non-consistency of the histogram density estimate from ergodic observations (cf. p.60. in \cite{GyHaSaVi89}).  
The dynamical system $(\Omega,\F,\mu,T)$  is determined by
 $\Omega=[0,1)$, 
$\F$ the Borel $\sigma$-algebra, $\mu$ the Lebesgue measure on $[0,1)$, and $T \omega=\omega+\alpha$ mod $1$ for some  irrational $\alpha$. 
 The dynamical system $(\Omega,\F,\mu,T)$ is stationary and ergodic by \cite{Billingsley}. Let $X_i(\omega)=T^{i+1}\omega$. 
We will apply Rohlin's lemma (cf. \cite{Shie73}), according to which 
if $(\Omega, \F,\mu ,T)$ is a nonatomic stationary and  ergodic dynamical system  then given 
$\epsilon>0$, and positive integer $N$, there exists a set $S\in \F$ 
such that 
$$
S, T^{-1} S, \dots, T^{-N+1} S
$$
are disjoint and 
$$
\mu(\cup_{i=0}^{N-1}T^{-i} S)\ge 1-\epsilon.
$$ 
For   $N=4n$ and $\epsilon=0.5$
we are assured of the existence of  a set $S\in \F$, such that 
$$
\mu(\cup_{i=0}^{4n-1} T^{-i} S)\ge 0.5.
$$ 
Put
$$B_n=\cup_{i=0}^{n-1} T^{-i}S$$
and
$$C_n=\cup_{i=0}^{2n-1}T^{-i}S.$$
Since $T^{-i} S$ $i=0,\dots,4n-1$ are disjoint and $T$ is measure preserving, we have 
$\mu(B_n)\ge 1/8$ and $1/4\le \mu(C_n)\le 1/2$. 
 Let  $X_i(\omega)=T^{i+1}\omega$.
The definitions of $B_n$ and $C_n$ imply that  all of 
$T^{-i} B_n\subset C_n$ 
for $i=0,\dots,n-1$ and  thus on the event $B_n$
all of the random variables $X_{-1},\dots,X_{-n}$ are in $C_n$, thus 
${1\over n} \sum_{i=1}^n I_{[X_{-i}\in C_n]}=1$.  
Now let $\{H_{n,j}\}_{j=1}^{q(n)}$ be a partition of the unit interval  by intervals of length $h_n=1/q(n)$ satisfying (\ref{condhntozero}) and (\ref{condnhntoinfty}). 
Let $A_{n,j}^+=H_{n,j}\cap C_n$ and 
$A_{n,j}^-=H_{n,j}\cap  \bar{C_n}$. 
Now let $\P_n=\{A_{n,j}^+,A_{n,j}^- \mbox{\  : \ $j=1,\dots,q(n)$}\}$. It is immediate that $\P_n$ satisfies conditions (\ref{r1}) and (\ref{r2}). 
\begin{eqnarray}
\nonumber
\lefteqn{
\int |\hat m_n(x)-m(x)|\mu(dx)}\\
\nonumber
&=& 
\sum_{j=1}^{q(n)} \int_{A_{n,j}^+} | {\nu_n(A_{n,j}^+)\over 
\mu_n(A_{n,j}^+)} -m(x)|\mu(dx) 
+  \sum_{j=1}^{q(n)} \int_{A_{n,j}^-} | {\nu_n(A_{n,j}^-)\over 
\mu_n(A_{n,j}^-)} -m(x)|\mu(dx) \\
\nonumber
&\ge& \sum_{j=1}^{q(n)} \int_{A_{n,j}^-} | {\nu_n(A_{n,j}^-)\over 
\mu_n(A_{n,j}^-)} -m(x)|\mu(dx) \\
&\ge& \sum_{j=1}^{q(n)} | \nu_n(A_{n,j}^-){ \mu(A_{n,j}^-) \over 
\mu_n(A_{n,j}^-)} -\int_{A_{n,j}^-} m(x)\mu(dx)|.  
\end{eqnarray}
On the event $B_n$, $\mu_n(\bar{C_n})=0$
and   consequently  $\mu_n(A_{n,j}^-)=\nu_n(A_{n,j}^-)=0$.  
Therefore on the event $B_n$, 
\begin{eqnarray}
\nonumber
\lefteqn{ \int |\hat m_n(x)-m(x)|\mu(dx)}\\
\nonumber
&\ge& \sum_{j=1}^{q(n)} \int_{A_{n,j}^-} m(x) \mu(dx) \\
\nonumber
&\ge& \sum_{j=1}^{q(n)} \mu(A_{n,j}^-) \inf_{x\in H_{n,j}}  
( (x+\alpha ) \mbox{  \ mod \ $1$}). 
\end{eqnarray}
For $1\le j\le q(n)$ let $g_n(j)=\inf_{x\in H_{n,j}}  
( (x+\alpha ) \mbox{  \ mod \ $1$})$ and $r_n(j)=\min\{l\ge 1: g_n(j)<l h_n\}$. Notice that function $r_n: \ \{1,\dots,q(n)\}\to \{1,\dots,q(n)\}$ is onto and invertible.  Since $\mu(\bar C_n)\ge 0.5$, on the event $B_n$, 
\begin{eqnarray}
\nonumber
\lefteqn{ \int |\hat m_n(x)-m(x)|\mu(dx)}\\
\nonumber
&\ge& \sum_{j=1}^{q(n)} \mu(A_{n,j}^-) g_n(j)\\
\nonumber
&\ge& \sum_{j=1}^{q(n)} \mu(A_{n,j}^-) (r_n(j)-1)h_n\\
\nonumber
&\ge&  \sum_{l=1}^{\lfloor \mu(\bar C_n)/h_n\rfloor } 
h_n (l-1)h_n\\
\nonumber
&\ge&  \sum_{l=1}^{\lfloor 0.5/h_n\rfloor } 
(l-1)(h_n)^2\\
\nonumber
&\ge&  0.5 (h_n)^2 ({1\over 2 h_n} -1)
({1\over 2 h_n} -2) \\
&\ge& {1\over 16} 
\end{eqnarray}
if $h_n< 1/12$. 
Since $h_n\to 0$, for large $n$,  on the event $B_n$, the $L_1$ error 
is at least $1/16$. That is, for large $n$,    
\[
P(\int |\hat m_n(x)-m(x)|\mu(dx)\ge 1/16)
\ge \mu(B_n)
\ge {1\over 8}.
\]
The proof of Theorem~\ref{Th4} is complete. 
\endpf

\noindent
{\bf Remark 3} Let process $\{X_n\}$ and the sequence of partitions $\{\P_n\}$ be as in 
Theorem \ref{Th4}. 
Set $Z_n=X_{n-1}$ and $Y_n=X_{n}+(1-\alpha)$ mod $1$. Define  $m(z)=E(Y_0|Z_0=z)$. 
It is easy to see that $m(z)=z$. Define  $\hat m_n(z)$ as in ~(\ref{defpartesthmnz}) 
with partition $\P_n$. 
The proof of Theorem~\ref{Th4} shows that the sequence of partitions $\{\P_n\}$ satisfies 
conditions (\ref{r1}) and (\ref{r2}) and   
\[
P(\int |\hat m_n(z)-m(z)|\mu(dz)\ge 1/16)\ge {1\over 8}. 
\]

\noindent
{\bf Acknowledgement} The authors wish to thank Paul Algoet for drawing their attention to Ryabko's paper \cite{Ryabko88}. The second athor thanks Benjamin Weiss for his comments. Comments from the referees have been extremely useful.


\begin{thebibliography}{99}

\bibitem{Algoet92}
P. H. Algoet, 
"Universal schemes for prediction, gambling and 
portfolio selection,"
{\sl Annals Probab.,} vol 20, pp. 901--941, 1992.
Correction: {\sl ibid.}, vol. 23, pp. 474--478, 1995.





\bibitem{Bailey76}
D. H. Bailey,
{\sl Sequential Schemes for Classifying and Predicting 
Ergodic Processes.} Ph. D. thesis, Stanford University, 1976.

\bibitem{Billingsley} P. Billingsley, {\sl Ergodic Theory  and Information.} Wiley, 1965. 

\bibitem{Cover75}
T. M. Cover,
"Open problems in information theory,"
in {\sl 1975 IEEE Joint Workshop on Information Theory,}
pp. 35--36. New York: IEEE Press, 1975.


\bibitem{DeGy83} L. Devroye and L. Gy\"orfi,
"Distribution-free exponential upper bound on the $L_1$ error of
partitioning estimates of a regression function", In {\sl
Proceedings of the Fourth Pannonian Symposium on Mathematical
Statistics}, Konecny F., Mogyor\'odi, J. and Wetz, W. Eds., pp.
67-76, Budapest, Akad\'emiai Kiad\'o, 1983.


\bibitem{Gyor91}
L. Gy\"orfi,
"Universal consistencies of regression estimate for unbounded regression functions,"
in {\sl Nonparametric functional estimation and related topics, ed. G. Roussas,}
pp. 329--338. Dordrecht: Kluwer Academic Publishers, 1991.

\bibitem{GyHaSaVi89} L. Gy\"orfi, Haerdle, W., Sarda, P., and Ph. Vieu, 
   {\em Nonparametric Curve Estimation from Time Series,} Springer Verlag,
   Berlin, 1989.

\bibitem{GyLu92} L.  Gy\"orfi and G. Lugosi, "Kernel density estimation
   from ergodic  sample is not universally consistent",  {\em Computational
   Statistics and Data Analysis}, {\bf 14}, pp. 437-442, 1992.

\bibitem{GyMa90} L. Gy\"orfi and E. Masry,  "The $L_1$ and $L_2$ strong
  consistency  of recursive kernel density estimation from time series", {\em
   IEEE Trans. on Information Theory}, {\bf 36}, pp. 531-539, 1990.


\bibitem{MoYaAl96}
G. Morvai, S. Yakowitz, and P. Algoet,
"Weakly convergent nonpapametric forecasting of stationary time series,"
  {\em IEEE Transactions on Information Theory}, {\bf 43}, pp. 483-498, 1997.

\bibitem{MoYaGy95}
G. Morvai, S. Yakowitz, and L. Gy\"orfi,
"Nonparametric inferences for ergodic, stationary time series,"
{\sl Annals of Statistics.}, vol. 24, pp. 370--379, 1996. 

\bibitem{Ornstein78}
D. S. Ornstein,
"Guessing the next output of a stationary process,"
{\sl Israel J. Math.,} vol. 30, pp. 292--296, 1978.

\bibitem{Ornstein74}
D. S. Ornstein,
 {\sl Ergodic Theory, Randomness, and Dynamical Systems.}
Yale University Press, 1974.

\bibitem{Rosen}
M. Rosenblatt, "Density estimates and Markov sequences," in
M. Puri, ed., {\sl Nonparametric Techniques in Statistical Inference,}
Cambridge University Press, Oxford, 1970. (pp. 199-210)


\bibitem{Rous}
G. Roussas, "Nonparametric estimation in Markov processes,"
{\sl Ann. Inst. Statist. Math.} vol. 21, pp. 73-87, 1969.

\bibitem{Ryabko88}
B. Ya. Ryabko,
"Prediction of random sequences and universal coding,"
{\sl Problems of Inform. Trans.,} 
vol. 24, pp. 87-96, Apr.-June 1988.

\bibitem{Shie73}
P.C. Shields, {\sl The Theory of Bernoulli Shifts}, The University of Chicago
Press, 1973.

\bibitem{Shields91}
P.C. Shields,
"Cutting and stacking: a method for constructing stationary processes,"
{\sl IEEE Transactions on Information Theory,}
vol. 37, pp. 1605--1614, 1991.

\bibitem{Wiener}
N. Wiener, {\sl Extrapolation, Interpolation and Smoothing of Stationary
Time Series}, the MIT Press, Cambridge, Mass., 1949.

\bibitem{ZiLe78}
J. Ziv and A. Lempel, "Compression of individual sequences by variable
rate coding," {\sl IEEE Trans. Inform. Theory,}
vol. IT-24, pp530-536, Sept. 1978.

\end{thebibliography}
\end{document}